\input eplain 
\magnification  1200
\ifx\eplain\undefined \input eplain \fi
\footnote{}{\bf{Typesetted on:  \timestamp}}
\baselineskip13pt 

\overfullrule=0pt


\font\smalsmalbf=cmbx8
\font\tenbi=cmmib10

\font\eightbi=cmmib9

\font\fivebi=cmmib5
\newfam\bmifam\textfont\bmifam=\tenbi\scriptfont
\bmifam=\eightbi\scriptscriptfont\bmifam=\fivebi

\def\bold #1{\fam\bmifam#1}
\def\expp{{\rm exp}} 
\font\smalltenrm=cmr8
\font\smallteni=cmmi8
\font\smalltensy=cmsy8
\font\smallsevrm=cmr6   \font\smallfivrm=cmr5
\font\smallsevi=cmmi6   \font\smallfivi=cmmi5
\font\smallsevsy=cmsy6  \font\smallfivsy=cmsy5
\font\smallsl=cmsl8      \font\smallit=cmti8

\def\smallfonts{\lineadj{80}\textfont0=\smalltenrm  \scriptfont0=\smallsevrm
                \scriptscriptfont0=\smallfivrm
    \textfont1=\smallteni  \scriptfont1=\smallsevi
                \scriptscriptfont0=\smallfivi
     \textfont2=\smalltensy  \scriptfont2=\smallsevsy
                \scriptscriptfont2=\smallfivsy
      \let\it\smallit\let\sl\smallsl\smalltenrm}

\font\eightbi=cmmib10 at 8pt

\font\smathbold=msbm5
\font\mathbold=msbm9 at 10pt
\font\smathbold=msbm8

\def\R{{\hbox{\mathbold\char82}}}
\def\sR{{\hbox{\smathbold\char82}}}

\def\e{\epsilon}

\def\d{{\hbox{\gothic\char78}}}


\def\imp{{\;\;\Longrightarrow\;\;}}
\def\lineadj#1{\normalbaselines\multiply\lineskip#1\divide\lineskip100
\multiply\baselineskip#1\divide\baselineskip100
\multiply\lineskiplimit#1\divide\lineskiplimit100}
\def\l{\ell}
\def\remark#1.{\medskip{\noin\bf Remark #1.\enspace}}
\def\endpf{$$\eqno/\!/\!/$$}

\def\pf#1.{\smallskip\noin{\bf  #1.\enspace}}

\def\noin{\noindent}
\def\D{\Delta}

\def\ds{\displaystyle}
\def\ts{\textstyle}

\def\e{\epsilon}\def\part{\partial_t}

\def\Rn{\R^n}

\def\ref#1{{\bf{[#1]}}}

\def\p{\partial}

\def\half{{\ts{1\over2}}}

\def \L{{\cal L}}

\def\W{W^{2,1}_\Delta}\def\Wo{W^{2,1}_{\Delta,0}}
\centerline{\bf  Optimal limiting embeddings for $ \Delta$-reduced Sobolev spaces in $\bold {L^1}$}

\bigskip \centerline{Luigi Fontana, Carlo Morpurgo}
\bigskip
\midinsert
{\smalsmalbf Abstract. }{\smallfonts We prove sharp embedding inequalities for certain reduced Sobolev spaces that arise naturally in the context of Dirichlet problems with $L^1$ data.  We  also find the optimal target spaces for such embeddings, which  in dimension 2  could be considered as limiting cases of the Hansson-Brezis-Wainger spaces, for the optimal embeddings of borderline Sobolev spaces $W_0^{k,n/k}$.

}

\endinsert\medskip\vskip3em
\centerline{\bf 1. Introduction}\bigskip

In this paper we are concerned with special kinds of the so-called {reduced Sobolev spaces}, namely the spaces defined by 
$$W_{\Delta}^{2,1}(\Omega)=\big\{u\in W_0^{1,1}(\Omega): \Delta u\in L^1(\Omega)\,\big\}\eqdef{@1} $$
and 
$$W_{\Delta,0}^{2,1}(\Omega)=\hbox{{\rm closure of }}\; C_0^\infty(\Omega)\;\;\hbox{{\rm in the norm}}\;\;\|\Delta u\|_1\eqdef{@2}$$
which we name {\it $\Delta-$reduced Sobolev spaces}. Here $\Omega$ is a bounded open set of $\R^n$, and $W_0^{k,p}(\Omega)$ denotes the closure of the set of $C^\infty$ functions compactly supported in $\Omega$, in the norm $\|u\|_{k,p}=\big(\sum_{|\alpha|\le k}\|D^\alpha u\|_p^p\big)^{1/p}$. The spaces $\W(\Omega)$ and $\Wo(\Omega)$ could be regarded as natural domains for the Dirichlet Laplacian, as an unbounded operator in $L^1(\Omega)$. Indeed, for $f\in L^1(\Omega)$ the problem $-\Delta u=f$ has a unique solution $u\in \W(\Omega)$, and if $\Omega$ is smooth enough then such $u$ is the limit of $C^\infty$ functions in $\Omega$ which are continuous up to the boundary, with 0 boundary value. The same considerations can be made for the $L^p$ versions $W^{2,p}_\Delta$ and $W^{2,p}_{\Delta,0}$, obtained by replacing $W_0^{1,1}$ with $W_0^{1,p}$ and $\|\Delta u\|_1$ with $\|\Delta u\|_p$ in \eqref{@1} and $\eqref{@2}$. There is, however, an important difference:  if $p>1$ then  $W^{2,p}_{\Delta,0}(\Omega)=W^{2,p}_0(\Omega)$, and (if $\Omega$ is smooth enough)  $W^{2,p}_\Delta(\Omega)=W^{1,p}_0(\Omega)\cap W^{2,p}(\Omega)$, but  these assertions are both false in the case $p=1$.  The reason for this is, essentially, that the $L^1$ norms of the second partial  derivatives cannot be controlled by $\|\Delta u\|_1$ (see for example [O]). For a general $L^1$ theory of second order elliptic equations see [BSt].

From the above discussion it should be apparent that such spaces are natural choices if one would like to study the summability properties of solutions of the Dirichlet problem in $L^1$, or the exceptional case $n=2$ of the Moser-Trudinger embedding $W^{2,{n\over2}}_0\hookrightarrow e^{L^{n\over n-2}}$.

In a recent paper [CRT] Cassani Ruf and Tarsi investigated sharp embedding properties of the $\Delta-$reduced spaces in (1) and (2), for smooth  $\Omega$. Among the
 main results of [CRT] are the sharp forms of  the embeddings
of $\W(\Omega)$ into the Zygmund space $L_\expp(\Omega)$, when $n=2$, and into the weak-$L^{n\over n-2}$ space $L^{{n\over n-2},\infty}(\Omega)$, when $n\ge 3$.
These spaces are defined by the quasi-norms 
$$\|u\|_{L_\expp}^*=\sup_{0<t\le |\Omega|}{u^*(t)\over1+\log{|\Omega|\over t}},\qquad\qquad \|u\|_{{n\over n-2},\infty}^*=\sup_{0<t\le |\Omega|} t^{n-2\over n} u^*(t),\eqdef{@0aa}$$
where $u^*$ denotes the decreasing rearrangement of $u$ on $(0,\infty)$, and the sharp forms of the embeddings derived in 
[CRT, Thm. 1, Thm. 2] are written as
$$\|u\|_{L_\expp}^*\le{1\over4\pi}\|\D u\|_1,\qquad n=2\eqdef{@0a}$$
$$\|u\|_{{n\over n-2},\infty}^*\le {1\over n^{n-2\over n}(n-2)\omega_{n-1}^{2/n}}\|\D u\|_1,\qquad n\ge 3\eqdef{@0b}$$
where $\omega_{n-1}$ denotes the volume of the $(n-1)$-dimensional unit sphere. The quantities on the left hand side are quasi-norms defining the spaces $L_\expp(\Omega) $ and $L^{{n\over n-2},\infty} (\Omega)$, respectively; the constants on the right hand sides are sharp, that is they cannot be replaced by smaller constants.

In [CRT] the following slightly better estimate is in fact obtained   for any $u\in\W(\Omega)$:
$$ u^*(t)\le N^*_{|\Omega|}(t)\|\D u\|_1,\qquad 0<t\le|\Omega|,\qquad n\ge 2\eqdef{n1}$$
where 
$$N_{|\Omega|}^*(t)=\cases{\ds{{1\over 4\pi}\log{|\Omega|\over t} }& if $n=2$\cr
{1\over n^{n-2\over n}(n-2)\omega_{n-1}^{2/n}}\big(t^{-{n-2\over n}}-|\Omega|^{-{n-2\over n}}\big)& if $n\ge3,$\cr} \eqdef{n2}$$
denotes the decreasing rearrangement of the Green function of the Laplacian for the ball of volume $|\Omega|$, with pole at the origin (see the proofs of Thm. 1,  Thm. 3 and Prop. 12 in [CRT]).

It must be noted that  inequality \eqref{n1} was obtained several years ago by Alberico and Ferone (see [AF] Theorem 4.1 and Remark 4.1 for the case $n=2$, and Theorem 5.1 for the case $n\ge 3$, which trivially yields \eqref{n1}). In [AF] it is in fact shown that $u^*(t)\le N_{|\Omega|}^*(t)\|P u\|_1$, for a general class of second order elliptic operators $P$, such that the Dirichlet problem  $Pu=f$
admits a unique weak solution $u\in L^1$, for each $f\in L^1(\Omega)$; in such generality, however, one cannot expect the inequality to be sharp. Related results are also contained in [AFT] and [Al].

Regarding the analogous results for the space $\Wo(\Omega)$, i.e.   the case of compactly supported functions,  only partial results were obtained in [CRT],  which however revealed an  intriguing aspect: among all functions of $\Wo(\Omega)$  which are either radial or nonnegative, inequalities \eqref{@0a}, \eqref{@0b}  and \eqref{n1} continue to hold but the constants are {\it halved}.  In particular,      Cassani, Ruf and Tarsi proved   that (see [CRT]  proofs of   Prop. 14 and Prop. 16) for any $t\in (0,|\Omega|\,]$ 
$$ u^*(t)\le {1\over 2} \, N^*_{|\Omega|}(t)\|\D u\|_1,\qquad u\in \Wo(\Omega),\quad u\ge0 {\hbox { or } } u {\hbox{ radial, }} \quad  n\ge 2\eqdef{n3}$$
and consequently ([CRT], Thm. 5, Prop. 14, Prop. 16)
$$\|u\|_{L_\expp}^*\le{1\over8\pi}\|\D u\|_1,\qquad n=2\eqdef{n4}$$
$$\|u\|_{{n\over n-2},\infty}^*\le {1\over 2n^{n-2\over n}(n-2)\omega_{n-1}^{2/n}}\|\D u\|_1,\qquad n\ge 3,\eqdef{n5}$$
for any $u\in\Wo(\Omega)$ which is either nonnegative or radial, and
with sharp constants,  within that class of functions.  One of the original motivations of this work was to find out whether the inequalities in  \eqref{n4}, and \eqref{n5} would still be valid, and therefore sharp, in the whole space $\Wo(\Omega)$. 

The first main result of this paper is the following sharp version of \eqref{n3}: if $\Omega$ is open and bounded, then for any $t\in(0,|\Omega|\,]$
$$ u^*(t)\le 2^{-2/n} N^*_{|\Omega|}(t)\|\D u\|_1,\qquad u\in\Wo(\Omega),\quad n\ge 2.\eqdef{n6}$$
 and the constant $2^{-2/n}$  in \eqref{n6} is sharp,  in the sense that it cannot be replaced with a smaller constant if $t$ is allowed to be sufficiently small. We will also prove sharpness of \eqref{n6}  for any given $t$ when $\Omega$ is either a ball ($n=2$) or the whole of $\R^n$ ($n\ge3$). As a consequence of \eqref{n6} we then find that allowing $u$ to be an arbitrary function in $\Wo(\Omega)$ (not just nonnegative or radial) inequality \eqref{n4} continues to hold, with sharp constant, whereas \eqref{n5} is replaced with 
$$\|u\|_{{n\over n-2},\infty}^*\le {2^{-2/n}\over n^{n-2\over n}(n-2)\omega_{n-1}^{2/n}}\|\D u\|_1,\qquad n\ge 3,$$
with  sharp constant. 

To prove \eqref{n6}, we will first rederive  \eqref{n3} (and also \eqref{n1}) for arbitrary open and bounded  $\Omega$,  as  a relatively straightforward consequence of Talenti's comparison theorem (which was also the starting point  in [AF] and [CRT])  and a well-known formula that goes back to Talenti [T],  for the solution of the Dirichlet problem on a ball with radial data (see \eqref{@15a}, \eqref{@16}).  The use of such formula in combination with Talenti's type comparison theorems  allows one to   obtain optimal norm estimates of  the solution $u$ of a Dirichlet problem $-\Delta u=~f$  in terms of norms of $f$; this idea  was already mentioned and used elsewhere (see for example [AFT, Prop. 3.1] and comments thereafter, and also [AF, proof of Theorem~4.1]).

The presence of the factor ${1\over2}$ in \eqref{n3} is perhaps better clarified in our proof, which is based on the simple observation that if $u$ is compactly supported in $\Omega$, then 
$\int_\Omega \Delta u=0$, and
$$\int_\Omega (\D u)^+ dx=\int_\Omega (\D u)^- dx={1\over 2}\,\|\D u\|_1\eqdef w$$
where $(\D u)^+$ and $(\D u)^-$ denote the positive and negative parts of $\D u$.  
 The proof of \eqref{n6}  will be then obtained by carefully combining estimates for the distribution functions of the positive and the negative parts of $u$. We will also introduce  natural families of radial extremal functions for \eqref{n1} and \eqref{n3}, essentially Green's potentials of normalized characteristic functions of balls or annuli; by suitably translating such functions  we will be able to produce a family of extremals for \eqref{n6}.

An immediate consequences of \eqref{n6} when $n=2$ is the following  Brezis-Merle type inequality   
$$\sup_{u\in\Wo(\Omega)}\int_{\Omega} e^{\alpha {|u(x)|\over\|\D u\|_1}}dx\le {8\pi\over 8\pi-\alpha}|\Omega|,\qquad \alpha<8\pi,\eqdef {bm1}$$
where the left-hand side is infinite if $\alpha=8\pi$, and  with sharpness of the constant ${8\pi\over 8\pi-\alpha}$ when $\Omega$ is a ball. The same  inequality holds for $\W(\Omega)$ with $8\pi$ replaced by $4\pi$: 
$$\sup_{u\in\W(\Omega)}\int_{\Omega} e^{\alpha {|u(x)|\over\|\D u\|_1}}dx\le {4\pi\over 4\pi-\alpha}|\Omega|,\qquad \alpha<4\pi,\eqdef{bm}$$ and as such it also appears in [AF, Thm 3.1], as a consequence of \eqref{n1}. The orignal Brezis-Merle inequality was obtained  in [BM] and it is essentially \eqref{bm}, but with a larger  right-hand side. Similar inequalities without explicit  right-hand side constants, but slightly more general integrands, where  also obtained in [CRT], but either on $\W(\Omega)$ or for functions of $\Wo(\Omega)$ which are nonnegative or radial. The Brezis-Merle inequality  quantifies the exponential integrability of functions  in $\Wo(\Omega)$ and $\W(\Omega)$, when $n=2$; indeed it is well known that the function $u$ is in $L_\expp(\Omega)$ if and only if $\int_\Omega e^{\lambda |u|} dx<\infty$, for some $\lambda>0$.

We observe that the discrepancy between the optimal ranges of $\alpha$'s in \eqref{bm1} and \eqref{bm} is a phenomenon that is peculiar to $L^1$ and the identities in \eqref w. Indeed, the analogous sharp exponential  inequality when $n>2$ 
$$\int_{\Omega}e^{\alpha \big({|u(x)|\over\|\D u\|_{n/2}}\big)^{n\over n-2}}dx\le C,\qquad 0<\alpha\le n(n-2)^{n\over n-2}\omega_{n-1}^{2\over n-2}\eqdef{mt}$$
was obtained by Adams [A] for the space $W_0^{2,n/2}(\Omega)=W_{\D,0}^{2,n/2}(\Omega)$, but it can be easily extended to the larger space $W_\Delta^{2,n/2}(\Omega)=W_0^{1,p}(\Omega)\cap W^{2,p}(\Omega)$, with the {\it same} sharp range of $\alpha$'s. The reason for that is that if $p>1$ then $\|f^+\|_p$ can be made arbitrarily close to $\|f\|_p$ within the class of functions with zero mean;  the vanishing of the mean of $\Delta u$ plays no role in \eqref{mt}, as opposed to the case $n=2$ in \eqref{bm1} and \eqref{bm},  where \eqref w causes a doubling of the largest exponential constant, going from
 general solutions of Dirichlet problems to compactly supported functions.

When $n\ge3$ one instead obtains, as a result of \eqref{n6}, an estimate of type 
$$\sup_{u\in \Wo(\Omega)}{\|u\|_q\over\|\D u\|_1}\le C(n,q,|\Omega|),\qquad 1\le q<{n\over n-2} $$
and a similar estimate for $\W(\Omega)$, using \eqref{n1}. In Corollary 2 we will exhibit a specific constant $C(n,q,|\Omega|)$ which is sharp in the case of $\W(\Omega)$ and  $\Omega$  a ball. Similar estimates without explicit constants, but slightly more general otherwise,  were also obtained in [CRT], but again, only on $\W(\Omega)$ or for functions of $\Wo(\Omega)$ which are nonnegative or radial.

A question of interest that one can raise, in view of the embedding results in [CRT] and in the present paper, is the following: {\it What is  the  smallest target space for the embeddings of $\Wo(\Omega)$?}

 A natural request in this sort of questions is that our admissible  target spaces be the so called {\it rearrangement invariant} spaces; those are Banach spaces $(X,\|\cdot \|_X)$ of Lebesgue measurable   functions on $\Omega$ with the property that $\|u\|_X=\|w\|_X$, whenever $u$ and $w$ are equimeasurable. This problem has been fully investigated in the case of   the classical Sobolev spaces embeddings. In particular, for the borderline embeddings of $W_0^{k,n/k}(\Omega)$, ($n>k$), the optimal r.i. target spaces turns out to be the so called Hansson-Brezis-Wainger spaces ([BW], [CP], [EKP], [H], [MP]); such spaces  are stricly contained in the exponential classes involved in  the Adams-Moser-Trudinger inequalities [A]. See also [C], where  optimal embedding results are obtained for general Orlicz-Sobolev spaces, including those of Hansson-Brezis-Wainger as special cases.  

The second  main result of this paper is that the optimal target space for the embedding $\Wo(\Omega)\hookrightarrow~X$, where $X$ is a r.i. space over $\Omega$, is the space of functions  $$L_{\expp,0}(\Omega)=\Big\{u\in L_\expp(\Omega): \;\lim_{t\to 0}{u^{**}(t)\over \log{1\over t}}=0\Big\},\qquad {\hbox{ when 
 }} \;n=2,$$and
$$L_0^{{n\over n-2},\infty}(\Omega)=\Big\{u\in L^{{n\over n-2},\infty}(\Omega): \;\lim_{t\to 0}t^{n-2\over n}u^{**}(t)=0\Big\},\qquad {\hbox { when }}\; n\ge 3,$$
where $u^{**}(t)={1\over t}\int_0^t u^*(s)ds$ denotes the so-called maximal function of $u^*$. It is easy to see that the limit conditions in the above spaces can be unified as
$$\lim_{t\to 0}{u^*(t)\over N_{|\Omega|}^*(t)}=0\eqdef {@0h}$$
which is obviously a stronger condition than \eqref{n6} from the point of view of ``best target space".
 
When $n=2$  the space $L_{\expp,0}(\Omega)$ is a Banach subspace of $L_\expp(\Omega)$, endowed with the norm
$$\|u\|_{L_\expp}=\sup_{0<t\le |\Omega|}{u^{**}(t)\over 1+\log {|\Omega|\over t}} \eqdef{@0j}$$
and our optimal embedding result can be interpreted as the limiting case of the optimal borderline embeddings obtained by  Hansson and Brezis-Wainger for $W_0^{k,n/k}(\Omega),\; n>k$.

When $n\ge 3$ the space $L_0^{{n\over n-2},\infty}(\Omega)$ is  a Banach subspace of $L^{{n\over n-2},\infty}(\Omega)$,  endowed with the norm
$$\|u\|_{{n\over n-2},\infty}=\sup_{0<t\le|\Omega|} t^{n-2\over n}u^{**}(t).\eqdef{@0k}$$

The space $L_{\expp,0}(\Omega)$  can also be characterized as the closure of the class of simple measurable functions on $\Omega$, in the norm $\|\cdot \|_{L_\expp}$, and also as the subspace of all order continuous elements of $L_\expp(\Omega)$ (i.e. those $f\in L_\expp(\Omega)$ such that if $|f_n|\le |f|$ and $|f_n|\downarrow 0$ then $\|f_n\|_{L_\expp}\downarrow 0$). This is also true for $L_0^{{n\over n-2},\infty}(\Omega)$, and in fact for any Marcinkievicz space $M_w(\Omega)$, defined by the norm $\|u\|_{M_w}=\sup\{u^{**}(t) w(t)\}$, for a quasiconcave function~$w$, and its subspace $M_w^0(\Omega)=\{u\in M_w(\Omega):\;\lim_{t\to0}u^{**}(t)w(t)=0\}$ (see for example [KPS], and also [KL] which contains  a nice summary of the properties of $M_w^0$).

It is important to note  that our optimal spaces $L_{\expp,0}$ and $L_0^{{n\over n-2},\infty}$  do not satisfy the so-called Fatou property, that is, they are not closed under a.e. limits of uniformly bounded   sequences. For this reason the definition of r.i. space that we adopt here, given for example in [KPS], is the more general one, which  does not require the Fatou property. It is an easy consequence of our result, however,   that the optimal r.i. spaces  {\it with} the Fatou property that contain $\Wo(\Omega)$ are $L_\expp(\Omega)$, when $n=2$, and $L^{{n\over n-2},\infty}(\Omega)$, when $n\ge3$ (see Theorem~2).

Our optimality results improve those obtained in Alberico-Cianchi, [AC], namely Theorem 1.1 in case $k=+\infty,\,n>p=2$ and Theorem 1.2, (iii), $k=+\infty,\;n=p=2$. In such theorems the authors prove in particular the optimality of the norms $\|u\|_{{n\over n-2},\infty}$ ($n\ge3$) and $\|u\|_{L_\expp}$ ($n=2$) in the inequality
 $$\|u\|_X\le C\|f\|_1\eqdef d$$
among all r.i. spaces $X$ satifying the Fatou property, assuming that the inequality is valid for all $f\in L^1$ and all solutions $u$  of a general class of boundary vaue problems, which includes the Dirichlet problem.  Their proof is based on a duality argument and the fact that if $X$ is an r.i. space with the Fatou property then its second associate space $X''$ coincides with $X$. It is well known that if $X$ does not satisfy the Fatou property, then $X$ is a proper subspace of $X''$ (see for example [BS], [KPS] and [KL] for a  summary of these and more facts on r.i. spaces, and references therein). In our result we assume only the minimal set of axioms for an r.i space, and the validity of \eqref d when $f=-\D u$, and $u$ compactly supported in $\Omega$, i.e. when $u\in \Wo(\Omega)$.

 Our proof is self-contained and borrows some ideas used in [CP, Thm. 5], for the spaces $W^{k,n/k}$. The key step is to prove that  for a function $u$ satisfying \eqref{@0h} and with support inside a ball of volume $V$ one has 
$$ u^{*}(t)\le (Tf)^*(t),\qquad 0<t\le V$$
where $T$ is the Green potential for the ball, and $f$ is a suitable  positive radial function on the ball. This is a version of [CP, Thm. 4] that is suited to our situation.
\bigskip
\centerline{\bf 2. Sharp embedding inequalities for $\W(\Omega)$ and $\Wo(\Omega)$}\bigskip

If $\Omega$ is an open set of $\R^n$ and  $u:\Omega\to\R$ is Lebesgue measurable, the decreasing rearrangement of $u$ is the function
$$u^*(t)=\inf\big\{s\ge0:\,|\{x\in \Omega:\,|u(x)|>s\}|\le t\big\},\qquad t>0$$
that is the function on $[0,+\infty)$ that is equimeasurable with $u$ and also decreasing.

On a ball $B_R=B(0,R)$ let 
$$N_{B_R}(r)=\cases{c_n (r^{2-n}-R^{2-n}) & if $n\ge3$\cr {\ds{1\over 2\pi}\log{R\over r}} & if $n=2$ \cr},\qquad 0<r\le R$$
with $$c_n={1\over (n-2)\omega_{n-1}},\qquad \omega_{n-1}={2\pi^{n/2}\over\Gamma\big({n\over2}\big)}$$

If $B_R$ is a ball of given volume $V$ and $0<t\le V$, we let 
$$N_V^*(t)=N_{B_R}\bigg(\Big({nt\over\omega_{n-1}}\Big)^{1/n}\bigg)=\cases{\ds{{1\over 4\pi}\log{V\over t} }& if $n=2$\cr
c_n\Big(\ds{\omega_{n-1}\over n}\Big)^{n-2\over n}\big(t^{-{n-2\over n}}-V^{-{n-2\over n}}\big)& if $n\ge3.$\cr}$$
 Note that if $G_{B_R}(x,y)$ is the Green function for the ball of volume $V$ then $N_V^*(t)$ is  the decreasing rearrangement of $G_{B_R}(x,0)$.

When $n\ge 3$ we also set
$$N_\infty^*(t)=c_n\Big(\ds{\omega_{n-1}\over n}\Big)^{n-2\over n}\,t^{-{n-2\over n}},\qquad t>0.\eqdef{@9w}$$

The $\Delta-$reduced spaces $\W(\Omega)$ and $\Wo(\Omega)$ are defined in \eqref{@1} and \eqref{@2}. Note that those definitions make sense for arbitrary open sets, not necessarily bounded. In particular when $\Omega=\R^n$ it is straightforward to check that $\W(\R^n)=\Wo(\R^n)$.
\medskip
\eject
\proclaim Theorem 1. Let $\Omega\subseteq \R^n$, $n\ge 2$, be open and bounded with volume $|\Omega|$. Then,\smallskip
\item{a)} For all $u\in \W(\Omega)$
$$u^*(t)\le N_{|\Omega|}^*(t)\|\D u\|_1,\qquad  0<t\le|\Omega|\eqdef{@10}$$
\item{b)} For all $u\in \Wo(\Omega)$ and $n\ge2$
$$u^*(t)\le {2^{-2/n}}\,N_{|\Omega|}^*(t)\|\D u\|_1,\qquad  0<t\le|\Omega|\eqdef{@11}$$
and if $n\ge 3$ and either $u\ge0$ or $u$ radial and $\Omega$ a ball,  then
$$u^*(t)\le {1\over 2}\,N_{|\Omega|}^*(t)\|\D u\|_1,\qquad  0<t\le|\Omega|\eqdef{@12}$$
When $n\ge 3$ both \eqref{@11} and \eqref{@12} hold for  $\Omega$ unbounded, with the convention in \eqref{@9w}.\smallskip
\item{c)} The inequalities in a) and b) are sharp in the following sense:
$$\sup_{u\in X,\,0<t\le |\Omega|}{u^*(t)\over N_{|\Omega|}^*(t)\|\D u\|_1}=\cases{1 & if $X=\W(\Omega)$\hskip11.2em{\rm \eqdef{@12a}}\cr \cr 2^{-2/n} & if $X=\Wo(\Omega)$\hskip11.05em{\rm \eqdef{@12b}}\cr\cr \ds{1\over2}& if  $X=\Wo(\Omega)\cap\{u {\hbox{  radial}}\}$\hskip5.5em{\rm \eqdef{@12c}}\cr &or $X=\Wo(\Omega)\cap\{u\ge0\}.$\cr}$$
Moreover, if  $B$ is any ball, then for each $t\in(0,|B|\,]$ 
$$\sup_{u\in X}{u^*(t)\over\|\D u\|_1}=\cases{N_{|B|}^*(t) & if $X=\W(B)$\hskip13.8em{\rm \eqdef{@12d}}\cr \cr  \ds{1\over2}N_{|B|}^*(t)& if  $X=\Wo(B)\cap\{u {\hbox{  radial}}\}$\hskip8.25em{\rm \eqdef{@12f}}\cr & or $X=\Wo(B)\cap\{u\ge0\}$\cr}$$
and also
$$\sup_{u\in \Wo(\sR^n)}{u^*(t)\over\|\D u\|_1}= 2^{-2/n}N_\infty^*(t). \eqdef{@12e} $$
\par
\medskip
\noin{\bf Remark.} As we noted in the introduction,  \eqref{@10}  appears in [AF] and  [CRT] and \eqref{@12} appears in [CRT], in case $\Omega$ is smooth. 
\smallskip
As an immediate consequence of Theorem 1 we obtain sharp norm embeddings for the spaces $\W(\Omega)$ and $\Wo(\Omega)$. Recall that 
$$L_\expp(\Omega)=\big\{u:\Omega\to\R,\;u {\hbox{  measurable  and }}\|u\|_{L_\expp}^*<\infty\big\}\eqdef{@G}$$
and
$$L^{{n\over n-2},\infty}(\Omega)=\big\{u:\Omega\to\R,\;u {\hbox{  measurable  and }}\|u\|_{{n\over n-2},\infty}^*<\infty\big\},\eqdef{@H}$$
where the quasi-norms $\|u\|_{L_\expp}^*$ and $\|u\|_{{n\over n-2},\infty}^*$ are defined as in \eqref{@0aa}. Note that in \eqref{@G}, \eqref{@H} the norms  $\|u\|_{L_\expp}$ and $\|u\|_{{n\over n-2},\infty}$ defined in \eqref{@0j}, \eqref{@0k} can be equivalently used in place of the corresponding quasi-norms.
\proclaim Corollary 1. Let $\Omega\subseteq \R^n$, $n\ge 2$, be open and bounded.
If $n=2$ then $\W(\Omega)\hookrightarrow L_\expp(\Omega)$ and  in particular
$$\|u\|_{L_{\expp}}^{*}\le {1\over 4\pi}\|\Delta u\|_1,\qquad w\in \W(\Omega)\eqdef{@x1}$$
$$ \|u\|_{L_{\expp}}^{*}\le {1\over 8\pi}\|\Delta u\|_1,\qquad w\in \Wo(\Omega)\eqdef{@x2}$$
and the constants ${1\over 4\pi}$ and ${1\over 8\pi}$ are sharp, i.e. they cannot be replaced by smaller constants.
\smallskip
If $n\ge 3$ then $\W(\Omega)\hookrightarrow L^{{n\over n-2},\infty}(\Omega)$ and in particular
$$\|u\|_{{n\over n-2},\infty}^{*}\le c_n\Big(\ds{\omega_{n-1}\over n}\Big)^{n-2\over n}\|\Delta u\|_1,\qquad w\in \W(\Omega)\eqdef{@x3}$$
$$ \|u\|_{{n\over n-2},\infty}^{*}\le 2^{-2/n}c_n\Big(\ds{\omega_{n-1}\over n}\Big)^{n-2\over n}\|\Delta u\|_1,\qquad w\in \Wo(\Omega)\eqdef{@x4}$$
and the constants are sharp.\par
\smallskip
\noin{\bf Remark.}  Corollary 1 continues to hold if $\|u\|_{L_\expp}^*$ and $\|u\|_{{n\over n-2},\infty}^{*}$ are replaced by the larger quantities $\|u\|_{L_\expp}$, $\|u\|_{{n\over n-2},\infty}$, and the constants in \eqref{@x1}-\eqref{@x4} are multiplied by ${n\over2}$. The reason for this is that 
 $$N_V^{**}(t)={1\over t}\int_0^t N_V^*(u)du=\cases{\ds{{1\over 4\pi}\Big(1+\log{V\over t} \Big)}& if $n=2$\cr
c_n\Big(\ds{\omega_{n-1}\over n}\Big)^{n-2\over n}\big({n\over2}\,t^{-{n-2\over n}}-V^{-{n-2\over n}}\big)& if $n\ge3$,\cr}$$
so that 
$N_V^{**}(t)\sim{n\over 2}N_V^*(t)$, as $t\to 0$.\smallskip
Another immediate consequence of the estimates of Theorem 1 are the following sharp versions of the Brezis-Merle and Maz'ya's  inequalities:\smallskip
\proclaim Corollary 2. Let $\Omega\subseteq \R^n$, $n\ge 2$, be open and bounded. If  $n=2$ then 
$$\int_{\Omega}e^{\alpha\,{|u(x)|\over \|\D u\|_1}}dx\le {4\pi \over 4\pi-\alpha}|\Omega|,\qquad 0<\alpha<4\pi,\;u\in \W(\Omega)\eqdef{@y1}$$
$$\int_{\Omega}e^{\alpha\,{|u(x)|\over \|\D u\|_1}}dx\le {8\pi \over 8\pi-\alpha}|\Omega|,\qquad 0<\alpha<8\pi,\;u\in \Wo(\Omega)\eqdef{@y2}$$
and  the integrals are infinite if $\alpha=4\pi$ in \eqref{@y1} and $\alpha=8\pi$ in \eqref{@y2}. If $\Omega$ is a ball the constants ${4\pi\over 4\pi-\alpha},\;{8\pi\over 8\pi-\alpha}$ are sharp.\smallskip

If $n\ge 3$ then, for $1\le q<{n\over n-2}$
$$\|u\|_q\le c_n\Big({\omega_{n-1}\over n}\Big)^{n-2\over n}\bigg[{\Gamma\big({n\over n-2}-q\big)\Gamma(q+1)\over\Gamma\big({n\over n-2}\big)}\bigg]^{1/q}\,|\Omega|^{{1\over q}-{n-2\over n}}\|\D u \|_1,\qquad \; u\in \W(\Omega)\eqdef{@y3}$$
$$\|u\|_q\le 2^{-{2\over qn}}c_n\Big({\omega_{n-1}\over n}\Big)^{n-2\over n}\bigg[{\Gamma\big({n\over n-2}-q\big)\Gamma(q+1)\over\Gamma\big({n\over n-2}\big)}\bigg]^{1/q}\,|\Omega|^{{1\over q}-{n-2\over n}}\|\D u \|_1,\qquad \; u\in \Wo(\Omega)\eqdef{@y4}$$
and if  $\Omega$ is a ball the constant is sharp in \eqref{@y3}.\par

\bigskip
\pf Proof of Theorem 1. The first step in the  proof of \eqref{@10} and \eqref{@12} is Talenti's comparison theorem, as in [AF] and [CRT], 
 and the following well known formula for the solution of the Dirichlet problem $-\Delta v= f$ on the ball $B_R$ and with radial data $f\in L^1(B_R)$:
$$v(|x|)=N_{B_R}(|x|)\int_{|y|\le |x|} f(y)dy+\int_{|x|\le |y|\le R} N_{B_R}(|y|)f(y)dy\eqdef{@15a}$$
or, in polar coordinates,
$$\eqalign{v(\rho)&=\omega_{n-1}N_{B_R}(\rho)\int_0^\rho f(r) r^{n-1}\,dr+\omega_{n-1}\int_\rho^R N_{B_R}(r)f(r)r^{n-1}\, dr\cr&=-\omega_{n-1}\int_\rho^R N_{B_R}'(r)dr\int_0^r f(\xi)\xi^{n-1}d\xi.\cr} \eqdef{@16}$$
Note that if either $f\ge0$ or $f$ decreasing with mean zero,  then $v(\rho)$ given  as in \eqref{@16} is decreasing.\smallskip
What we need here is the following version of Talenti's result: let $\Omega$ be open and bounded and let $f\in L^1(\Omega)$  and let $f^\sharp(x)=f^*(|B_1||x|^n)$, the Schwarz symmetrization of $f$, supported in the ball $B_R$ with volume $|\Omega|$; if $u,v\in W_0^{1,1}(\Omega)$ are the unique solutions of $-\D u=f$ and $-\Delta v=f^\sharp$, then $u^*(t)\le v^*(t)$ for $t>0$. This result (including existence and uniqueness of the solutions)  follows by a routine  argument: 1) approximate  $f$ in $L^1$ via a sequence of $f_n\in C_0^\infty(\Omega)$; 2) solve the problems $-\D u_n=f_n$ , $-\D v_n=f_n^\sharp$; 3) use the uniform gradient estimate $\|\nabla u_n\|_1\le \|\nabla v_n\|_1\le C\|f\|_1$ (the left inequality for example is in [T, p. 715]) ; 4) show that $\{u_n\}$ is a Cauchy sequence convergent to $u$, the solution of $-\D u=f$; 5) apply Talenti's classical result to the $u_n$, and pass to the limit.
  
To prove \eqref{@10} we then  apply the above version of  Talenti's  theorem to a function $\in \W(\Omega)$, and conclude that $u^*(t)\le v^*(t)$ for $t>0$, where $v$ is the solution of $-\Delta v=(\Delta u)^\sharp$, $v=0$ on $\p B_R$. Next, note that the solution of $-\Delta v=f$ ($v=0$ on $\p B_R$) with $f$ radial given in \eqref{@15a} satisfies
$$|v(|x|)|\le N_{B_R}(|x|)\|f\|_1$$
which instantly gives \eqref{@10}.
\eject
A small modification of the above argument yields \eqref{@12} in the case $u\in\Wo(\Omega)$ with either $u\ge0$ or $u$ radial. Indeed,  assuming WLOG that $u\in C_0^\infty(\Omega)$, then $\int_\Omega \Delta u=0$, so letting $f=-\Delta u$, and $f^+, f^-$ the positive and negative parts of $f$, we have $\int_\Omega f^+=\int_\Omega f^-=\half\|f\|_1$. 
If $u$ is radial then \eqref{@15a} yields
$$-N_{B_R}(|x|)\int_{B_R}f^-(y)dy \le v(|x|)\le N_{B_R}(|x|)\int_{B_R} f^+(y)dy$$
or 
$$|v(|x|)|\le \half N_{B_R}(|x|)\|f\|_1$$
from which \eqref{@12} follows. If $u\ge0$ then letting $w$ be the solution of $-\D w=f^+$ on $\Omega$, with $w\in W^{1,1}_0(\Omega)$ we have $0\le u\le w$, by the maximum principle, and the result follows from part a) applied to $w$.

To prove \eqref{@11} we argue as follows. First, note that it is enough to prove the result for $u\in C_0^\infty(\Omega)$. For such given $u$  and  for each $\e\ge0$ consider the open subsets of $\Omega$  
$$\Omega_{\e}=\{x\in\Omega: u(x)>\e\},\qquad\Omega_\e'=\{x\in\Omega: -u(x)>\e\}$$
and the functions 
$$u_\e:=(u-\e)\big|_{\Omega_\e},\qquad u_\e'=(-u-\e)\big|_{\Omega_\e'}.$$

Sard's theorem combined with the implicit function theorem guarantee that for a.e. $\e>0$ both $\partial\Omega_\e$ and $\partial \Omega_\e'$ are smooth $C^\infty$ $(n-1)-$dimensional  manifolds; therefore, for each    such $\e$  both $u_\e$ and $u_\e'$ are $C^\infty$ in their domains, continuous up to the boundaries, and with zero boundary values, and if $f=-\D u$ they clearly solve the Dirichlet problems $-\Delta u_\e=f$ and $-\D u_\e'=-f$ in their domains. Let now  $w_e,\, w_\e'$ be the solutions to the Dirichlet problems
$$\cases{-\D w_\e=f^+ & on $\Omega_\e$ \cr w_\e=0   & on $\p \Omega_\e$\cr}\qquad \cases{-\D w_\e'=f^- & on $\Omega_\e'$ \cr w_\e'=0   & on $\p \Omega_\e'.$\cr}$$
Then we have $0\le u_\e\le w_\e$ and $0\le u_\e'\le w_\e'$, and 
also $w_\e\in W_\Delta^{1,2}(\Omega_\e)$, $w_\e'\in W_\Delta^{1,2}(\Omega_\e')$. We can then apply part a)
to deduce
$$(u_\e)^*(t)\le (w_\e)^*(t)\le N_{|\Omega_\e|}^*(t) \int_{\Omega_\e} f^+ dx, $$
for $0<t\le |\Omega_\e|$
and hence for $0<t\le |\Omega_0|$.
All the quantities involved above are monotone decreasing w.r. to $\e$ hence we deduce
$$(u_0)^*(t)\le N_{|\Omega_0|}^*(t)\int_{\Omega_0}f^+={1\over 2} N_{|\Omega_0|}^*(t)\|\D u\|_1\qquad 0<t\le |\Omega_0|\eqdef{@17}$$
Likewise, arguing with $u_\e',w_\e'$, we obtain  
$$(u_0')^*(t)\le {1\over 2} N_{|\Omega_0'|}^*(t)\|\D u\|_1\qquad 0<t\le |\Omega_0'|\eqdef{@18}$$
Let now $\lambda_V(s)$ be the distribution function of $N_V^*$, i.e. 
$$\lambda_V(s)=|\{t>0: N_V^*(t)>s\}|=\cases{V e^{-4\pi s} & if $n=2$ \cr \Big(\alpha_n s+V^{-{n-2\over n}}\Big)^{-{n\over n-2}} & if $n\ge 3$\cr}$$
where $\alpha_n=(n-2)n^{n-2\over n}\omega_n^{2/n}$.
With this notation we have, for $s>0$, 
$$\eqalign{|\{x\in\Omega:\;|u(x)|>s\}|&=|\{x\in \Omega_0:\; u_0(x)>s\}|+|\{x\in \Omega_0':\; u_0'(x)>s\}|\cr&\le
 \lambda_{|\Omega_0|}\Big({2s\over\|\D u\|_1}\Big)+\lambda_{|\Omega_0'|}\Big({2s\over\|\D u\|_1}\Big).\cr}\eqdef{@18}$$
Now note that $|\Omega_0|+|\Omega_0'|=|\Omega|$ and that
 $$\lambda_{|\Omega_0|}\Big({2s\over\|\D u\|_1}\Big)+\lambda_{|\Omega_0'|}\Big({2s\over\|\D u\|_1}\Big)\le 
\cases{|\Omega| e^{-8\pi s/\|\D u\|_1} & if $n=2$\cr
\Big( 2^{2/n}\ds{\alpha_n s\over\|\D u\|_1}+|\Omega|^{-{n-2\over n}}\Big)^{-{n\over n-2}} & if $n\ge 3$,\cr}\eqdef{@19}$$
since for $n=2$ there actually is equality, whereas for  $n\ge 3$ the right hand side of \eqref{@18} is maximized precisely when $|\Omega_0|=|\Omega_0'|=\half |\Omega|$. Inequalities \eqref{@18} and \eqref{@19} imply \eqref{@11}.

Now let us prove the sharpness statements. Introduce the radially decreasing functions
$$F_{\delta}^R={\chi_{B_\delta}^{}\over|B_\delta|},\quad 0<\delta<R$$
$$F_{\delta,\e}^R={\chi_{B_\delta}^{}\over2|B_\delta|}-{\chi_{A_{\e,R}}^{}\over2|A_{\e,R}|},\qquad 0<\delta<R-2\e<R$$
where
$$B_\delta=\{x:\, |x|\le \delta\},\qquad A_{\e,R}=\{x: R-2\e<|x|<R-\e\}.$$
Applying formula \eqref{@15a} we obtain that the solution $U_\delta^R$  of the Dirichlet problem
$$\cases{-\D U_{\delta}^R=F_{\delta}^R & on $B_R$\cr U_{\delta}^R=0 & on $\partial B_R$\cr}$$
is given by
$$U_\delta^R(x):=\cases{\ds{|x|^n\over \delta^n} N_{B_R}(|x|)+\ds{{1\over|B_\delta|}\int_{|x|<|y|<\delta}N_{B_R}(|y|)dy} & if $|x|<\delta$ \cr N_{B_R}(|x|) & if $\delta\le |x|\le R.$\cr}$$
which is nonnegative, radial and decreasing, so that 
$$(U_{\delta}^R)^*(t)=N_{|B_R|}^*(t),\qquad  |B_\delta|\le t\le |B_R|,$$ and this takes care of \eqref{@12d} immediately, since $U_\delta^R\in \W(B_R)$.

 If $\Omega$ is an arbitrary open and bounded set, then we can assume that $0\in\Omega$, and find $R$ so that $B_R\subseteq \Omega$. The function $U_\delta^R$ (extended to be 0 outside $B_R$) is not in $\W(\Omega)$, however we can argue that since $F_\delta^R\ge 0$ then the solution $U_\delta\in \W(\Omega)$ of $-\D U_\delta=F_\delta^R$ is nonnegative on $\Omega$  and satisfies $U_\delta^R\le U_\delta$ on $B_R$, by the maximum principle;  hence $(U_\delta)^*(t)\ge (U_\delta^R)^*(t)=N_{|B_R|}^*(t)$, for $|B_\delta|\le t\le |B_R|$.
It's then clear that taking $\delta_t$ so that $|B_{\delta_t}|=t$ gives
$${(U_{\delta_t})^*(t)\over N_{|\Omega|}^*(t)}\ge{N_{|B_R|}^*(t)\over N_{|\Omega|}^*(t)} \to 1,\qquad t\to 0,$$
thereby proving \eqref{@12a}.

Likewise,  the solution $U_{\delta,\e}^R$ to 
$$\cases{-\D U_{\delta,\e}^R=F_{\delta,\e}^R & on $B_R$\cr U_{\delta,e}^R=0 & on $\partial B_R$\cr}$$
can be computed explicitly, however all we need is that $U_{\delta,\e}^R$ is nonnegative, radial, decreasing on $(0,|B_R|]$, and 
$$ U_{\delta,\e}^R(x)=\cases{\half N_{B_R}(|x|)-\ds{1\over 2|A_{\e,R}|}\int_{A_{\e,R}} N_{B_R}(|y|)dy & if $\delta\le |x|\le R-2\e$\cr 0 & if $R-\e\le |x|\le R$\cr}\eqdef{@20}$$
all of which can be readily checked.  We then have $U_{\delta,\e}^R\in \Wo(B(0,R))$, and the above identity leads to \eqref{@12f}, since
$$\lim_{\e\to0}{1\over 2|A_{\e,R}|}\int_{A_{\e,R}} N_{B_R}(|y|)dy=0$$
For an arbitrary open and bounded $\Omega$, we can prove \eqref{@12c} like before, assuming $0\in \Omega$, $B(0,R)\subseteq\Omega$, this time observing that $U_{\delta,\e}^R\in \Wo(B(0,R))\subseteq \Wo(\Omega)$.
\def\d{\delta}
\def\xl{x_\lambda}
It remains  to settle \eqref{@12b} and \eqref{@12e} for $n\ge 3.$ We consider the functions
$$V_{\delta,\lambda}^R(x)=U_{\delta, R/4}^R(x)-U_{\delta,R/4}^R(x-\xl ),\qquad x_\lambda:=(\lambda,0,0...,0),$$
with
$$ \delta<\min\{\half,\half R\},\qquad\delta<\half\lambda<\half R,\eqdef{@12j}$$
so that \def\l{\lambda}
$$-\D V_{\delta,\l}^R={1\over 2|B_\d|}(\chi_{B_\d}^{}-\chi_{\xl+B_\delta}^{})-h_{\l}^R,$$
where $B_\d$ and $B_\d+\xl$ are disjoint  and where 
$$h_{\lambda}^R={1\over |A_{R/4,R}|}(\chi_{A_{R/4,R}}^{}-\chi_{\xl+A_{R/4,R}}^{})$$
which converges to 0 pointwise and in $L^1$, as $\l\to0$ for fixed $R$, and as $R\to+\infty$ for fixed $\l$; moreover, $|h_\d^R|\le C R^{-n}$ and
$$\int_{\sR^n}|h_\l^R|\le C{\l\over R}.\eqdef{@12k}$$
  Note that $V_{\d,\l}^R\in \Wo(B(0,R+\l))$.

In order to estimate the distribution function of $V_{\d,\l}^R$ on a given  $\Omega$ containing the support of such function,  write for $s>0$
$$|\{x\in\Omega: |V_{\d,\l}^R(x)|>s\}|\ge2|\{x:\,\delta<|x|<\half R,\,x_1<\half \l,\; |V_\d^R(x)|>s\}|$$
Note that \eqref{@20} gives 
$$U_{\d,R/4}^R(x)=\half c_n|x|^{2-n}-d_nR^{2-n},\qquad\d\le |x|\le \half R\eqdef{@21}$$
 for some $d_n>0$.

If    $x_1<\half\l$ then  $0\le U_{\d,R/4}^R(x-x_\l)\le U_{\d,R/4}^R(\half\xl)$, since $U_{\d,R/4}^R(x-x_\l)$ is radial decreasing about $x_\l$, and since $\d<\half \l<\half R$ we also have, using \eqref{@21}, 
$$|V_{\d,\l}^R(x)|\ge U_{\d,R/4}^R(x)-U_{\d,R/4}^R(\half x_\l) =\half c_n|x|^{2-n}-2^{n-3}c_n\l^{2-n},$$
and it is clear that the right-hand side is greater than $s$ if and only if $|x|<|x^*|$,
where 
$$|x^*|=\Big({2s\over c_n} +2^{n-2}\l^{2-n}\Big)^{-{1\over n-2}}<{\lambda\over2}<{R\over2}.$$
Conversely, if $|x^*|$ defined by the above equation satisfies $|x^*|\ge\d$, then 
$$|\{x\in \Omega:\,\half c_n|x|^{2-n}-2^{n-3}c_n\l^{2-n}>s\}|={\omega_{n-1}\over n} |x^*|^n={\omega_{n-1}\over n} \Big({2s\over c_n} +2^{n-2}\l^{2-n}\Big)^{-{n\over n-2}}.$$
Since $|x^*|\ge\d$ if and only if $s\le\half c_n\big( \d^{2-n}-2^{n-2}\l^{2-n}\big)>0$ (due to \eqref{@12j}), we finally obtain that for any such $s$
$$|\{x\in\Omega: |V_{\d,\l}^R(x)|>s\}|\ge2{\omega_{n-1}\over n}\Big[ \Big({2s\over c_n} +2^{n-2}\l^{2-n}\Big)^{-{n\over n-2}}-\delta^n\Big],\eqdef{@22}$$
which implies
$$(V_{\d,\l}^R)^*(t)\ge {c_n\over 2}\bigg[\bigg({nt\over 2\omega_{n-1}}+\delta^n\bigg)^{-{n-2\over n}} -2^{n-2}\l^{2-n}\bigg],\qquad 0\le t\le 2|B_{\lambda/2}|-2|B_\delta|.\eqdef{@23}$$

For a given open and bounded  $\Omega$, assume $0\in \Omega$, and fix $R<1$ so that $B(0,2R)\subseteq\Omega$. Pick any $\sigma$ with $0<\sigma<1/2$, and take $\delta<R^{1/\sigma}$  and $\lambda=\d^\sigma$, so that 
 $V_{\d,\d^\sigma}^R\in \Wo(B(0,2R))\subseteq\Wo(\Omega)$, and $\|\D V_{\d,\d^\sigma}^R\|_1\le 1+\|h_{\d^\sigma}^R\|_1\to 1$, as $\d\to0$. Therefore, \eqref{@23}  with $\d_t$ such that $2|B_{\delta_t}|=t^2$, and $t$ so small so that $\delta_t^\sigma>2\delta_t$, gives
$${(V_{\d_t\,\d_t^\sigma}^R)^*(t)\over N_{|\Omega|}^*(t) \|\D V_{\d_t,\d_t^\sigma}^R\|}\ge
2^{-2/n} c_n\bigg({ \omega_{n-1}\over n}\bigg)^{{n-2\over n}} {(t+t^2)^{-{n-2\over n}}-Ct^{-2\sigma{n-2\over n}}\over N_{|\Omega|}^*(t)(1+\|h_{\d_t^\sigma}^R\|_1)}\to2^{-2/n}$$
as $t\to0$, proving \eqref{@12b}.

If instead we fix $t>0$, then take $\Omega=\R^n$,   $\d_R$ so that $2|B_{\d_R}|=1/R$, and $R>1$ so large that  if $\l= R^\sigma$, with $0<\sigma<1$, then $t<2|B_{R^\sigma/2}|-2|B_{\d_R}|$, so that from \eqref{@12k} and \eqref{@23} we have 
$${(V_{\d_R, R^\sigma}^R)^*(t)\over \|\D V_{\d_R,R^\sigma}^R\|}\ge 2^{-2/n} c_n\bigg({ \omega_{n-1}\over n}\bigg)^{{n-2\over n}} {\bigg((t+R^{-1})^{-{n-2\over n}}-CR^{\sigma(2-n)}\bigg)\over (1+\|h_{R^\sigma}^R\|_1)}\to2^{-2/n}N_\infty ^*(t),$$
as $R\to+\infty$, yielding  \eqref{@12e}.
\endpf
\medskip
\pf Proofs of Corollaries 1,2. The inequalities \eqref{@x1}-\eqref{@y4} are straightforward consequences of  \eqref{@10} and \eqref{@11}. The proof of the sharpness statements can be easily obtained arguing as in the proof of Theorem 1, using the families of functions  $U_\d\in \W(\Omega),\,U_{\d,\e}^R\in \Wo(\Omega),\,V_{\d,\lambda}^R\in  \Wo(\Omega)$.\endpf
\noin{\bf Remark.} The question of the sharpness of \eqref{@y4} remains unsettled. The extremal families used in the above proofs seem to be unsuited for the computation of  the supremum of 
$|\Omega|^{-{1\over q}+{n-2\over n}}\|u\|_q\|\D u\|_1^{-1}$, over all open and bounded $\Omega$ and all $u\in\Wo(\Omega)$. 

\def\L{\Lambda}
\bigskip
\centerline{\bf 3. Optimal target spaces}
\bigskip
In this section we improve the embedding results of Corollary 1 from the point of view of ``smallest target space".
For $\Omega\subseteq\R^2$ define the space
$$L_{\expp,0}(\Omega)=\Big\{u\in L_\expp(\Omega): \;\lim_{t\to 0}{u^{**}(t)\over \log{1\over t}}=0\Big\}$$
which is a closed subspace  of $L_\expp(\Omega)$, endowed with the norm $\|u\|_{L_\expp}$.
Likewise, for $\Omega\subseteq \R^n$, $n\ge3$ define 
$$L_0^{{n\over n-2},\infty}(\Omega)=\Big\{u\in L^{{n\over n-2},\infty}(\Omega): \;\lim_{t\to 0}t^{n-2\over n}u^{**}(t)=0\Big\},$$
which is a closed subspace of $L^{{n\over n-2},\infty}(\Omega)$, endowed with the norm $\|u\|_{{n\over n-2},\infty}$.

Given a Lebesgue measurable set $\Omega$ let \def\M{{\cal M}}$\M_\Omega$ be the set of all Lebesgue measurable functions $f:\Omega\to[-\infty,\infty]$ which are a.e. finite (with the usual convention that a.e. equal functions are identified). A {\it rearrangement invariant} (r.i.) space over $\Omega$ is a Banach space $(X,\|\cdot\|_X)$ which is a subspace of $\M_\Omega$ satisfying the two properties
\smallskip
\noin (i) $\;|g|\le |f|$ a.e. and $f\in X$ $\imp  g\in X$ and $\|g\|_X\le \|f\|_X$ ($X$ is an ideal Banach lattice)
\smallskip
\noin{(ii)} if $f,g\in \M_\Omega$ are equimeasurable (i.e. if $|\{x\in \Omega: \;|f(x)|>s\}|=|\{x\in\Omega:\,|g(x)|>s\}|$ for each $s\ge0$),  then $\|f\|_X=\|g\|_X.$\smallskip

  In addition, we say that an r.i. space $(X,\|\cdot\|_X)$ satisfies the {\it Fatou property} if the following condition holds:
\smallskip
\noin {(iii)} if $0\le f_n\uparrow f$ a.e., with $f_n\in X$ and $\sup_n\|f_n\|_X<\infty$, then $f\in X$ and $\|f_n\|_X\uparrow \|f\|_X$.
 \smallskip
The Fatou property is easily seen to be equivalent to 

\noin{(iii')} if $f_n\to f$ a.e., with $f_n\in X$ and $\sup_n\|f_n\|_X<\infty$, then $f\in X$ and $\|f\|_X\le \liminf_n \|f_n\|_X$.

The above definition  of rearrangement invariant space is taken from [KPS] (where  it is  called ``symmetric space"); in other standard references, such as  [BS], the Fatou property is instead included in the defining axioms.
 
Clearly, both  $L_\expp(\Omega)$ and $L^{{n\over n-2},\infty}(\Omega)$ are 
rearrangement invariant spaces  over $\Omega$, and both of them satisfy the  Fatou property. The spaces $L_{\expp,0}(\Omega)$ and $L_0^{{n\over n-2},\infty}(\Omega)$ are  r.i. spaces over $\Omega$ which {\it do not} satisfy the Fatou property. This is easily seen by considering truncations of  the function $f(x)=N_B(|x|)$, where $B$ is any small ball inside $\Omega$. 

It is a known fact ([KPS, Thm 4.1])  that if conditions (i) and (ii) holds, $X$ is nontrivial,  and if $|\Omega|<\infty$ then
$$L^\infty(\Omega)\hookrightarrow X\hookrightarrow L^1(\Omega)$$
in the sense of continuous embeddings. The closed graph theorem also implies that any Banach space $Y$ which is a subset of  an r.i. space $X$ over $\Omega$, with $|\Omega|<
\infty$, is  continuously embedded in~$X$.
\medskip
\proclaim Theorem 2. For $n\ge2$ and $\Omega$ open and bounded in $\Rn$, let $\L_2^0(\Omega)=L_{\expp,0}(\Omega)$ if $n=2$, and $\L_n^0(\Omega)=L_0^{{n\over n-2},\infty}(\Omega)$ if $n\ge 3$. Then,  we have 
 $$\Wo(\Omega)\subseteq\W(\Omega)\subseteq \Lambda_n^0(\Omega)\eqdef{@14}$$
and for any rearrangement invariant  space $(X,\|\cdot\|_X)$ over $\Omega$
$$ \Wo(\Omega)\subseteq X\imp\Lambda_n^0(\Omega)\subseteq X.\eqdef{@15}$$
In other words, $\Lambda_n^0(\Omega)$ is the smallest target space $X$ for the embedding $\Wo(\Omega)\subseteq X$, among all r.i. spaces $X$.\smallskip
Moreover, if $(X,\|\cdot\|_X\|)$ is any r.i. space with the Fatou property (iii), then for $n=2$ 
$$ \Wo(\Omega)\subseteq X\imp L_\expp(\Omega)\subseteq X,\eqdef{@15z}$$
and for $n\ge 3$
$$ \Wo(\Omega)\subseteq X\imp L^{{n\over n-2},\infty}(\Omega)\subseteq X.\eqdef{@15zz}$$
\par

\smallskip\pf Proof of Theorem 2. If $u\in \W(\Omega)$, then the fact that $u\in \L_n^0(\Omega)$ follows easily from Talenti's comparison theorem combined with \eqref{@16}.

Let  now $(X,\|\cdot\|_X)$ be an r.i. space over  $\Omega$, endowed with the Lebesgue measure,  such  that $\Wo(\Omega)\subseteq X$. We claim that for any $u\in \L_n^0(\Omega)$ there exists a function $v\in \Wo(\Omega)$ and a constant $C$  such that 
$$u^{*}(t)\le v^*(t)+C,\qquad 0<t\le |\Omega|,\eqdef{bound}$$
which    implies  that $u\in X$ and therefore \eqref {@15}; obviously it is enough to show this  for $u\ge0$.

To prove the claim,  let us assume first WLOG that $0\in\Omega$ and that $u_0\in \L_n^0(\Omega)$ has support inside a ball $B_R\subseteq \Omega$. We now show that  we can find a nonegative integrable  function  $h:[0,|B_R|\,]\to\R$ such that 
$$N_{|B_R|}^*(t)\int_0^t h(s)ds\ge (u_0)^*(t),\qquad 0<t\le |B_R|.\eqdef{@24}$$
To prove the claim, let $g(t)=(u_0)^{**}(t)/N_{|B_R|}^*(t)$ ($0<t\le|B_R|)$, which is continuous and  converges to $0$ as $t\to0$ (by hypothesis), and let $f(t)=\sup_{0<s<t} g(s)$. This $f$ 
 is continuous, nonnegative, increasing, satisfies $f\ge g$, and $f(t)\to0$ as $t\to0$. 

Take any nonnegative, differentiable and decreasing function $m:(0,|B_R|\,]\to \R$, with $m(|B_R|)=0$ and $m(t)\to+\infty$ as $t\to0$ (for example $m(t)=\log(|B_R|/t)$), and let 
$$k(t)=-{1\over m(t)}\int_t^{|B_R|} f(s)m'(s)ds={1\over m(t)}\int_0^{m(t)} f(m^{-1}(u))du,\qquad 0<t<|B_R|;\eqdef k$$
such $k$ is differentiable, positive,  increasing, $k(t)\to f(|B_R|)$ if $t\to|B_R|$, and $k(t)\to0$ as $t\to0$. Therefore, the function  $h(t):=k'(t)$ is integrable, nonnegative  and it satisfies \eqref{@24}.

Now let us go back to our $u\in \L_n^0(\Omega)$, and assume that $u\ge0$, $u$ is not 0 a.e., $0\in \Omega$, and   $\lambda>0$ is such that
$|\{x\in\Omega:\,u(x)>\l\}|=|B_R|$, with $B_{2R}\subseteq \Omega$. Define 
$$u_0(x)=\max\{u(x),\l\}-\l=u(x)-\min\{u(x),\l\},\qquad x\in\Omega.$$

Clearly $(u_0)^*(t)=u^*(t)-\l$ for $0<t<|B_R|$, and $(u_0)^*(t)=0\ge u^*(t)-\l$ for $|B_R|\le t\le |\Omega|$, so that $u_0\in  \L_n^0(\Omega)$ and $u^*\le(u_0)^*+\lambda$. If $u_0^\#(x)=(u_0)^*(|B(0,x)|)$, for $x\in B_R$ and $u_0^\#(x)=0$ for $x\in\Omega\setminus B_R$, then  $u_0^\#\in \L_n^0(\Omega)$,  and $u_0^\#$ is supported in $B_R$.
Let  $f$ be the radial and integrable function on $B_R$ defined as  $f(x)=h(|B(0,x)|)$ with $h=k'$ and $k$ as in \eqref k. If  $v_0\in \W(B_R)$ is the solution of the problem $-\D v_0 = f$ given as in \eqref{@16}, then 
    by \eqref{@24}  
$$(v_0)^*(t)=N_{|B_R|}^*(t)\int_0^t h(s)ds+\int_t^{|B_R|}N_{|B_R|}^*(s) h(s)ds\ge u_0^*(t),
 \qquad 0<t\le|B_R|.\eqdef {k1}$$ On the other hand, if $v_1\in \W(B_{2R})$ solves $-\D  v_1=f$, (with $f=0$ outside $B_R$) with $v_1=0$ on $\p B_{R}$,  then $v_0\le v_1$ on $B_{R}$ (since $f\ge 0$), and we can construct a function $v\in \Wo(B_{2R})$, so that $v_1\le v+C$ for some constant $C$.
In order to do that, it is enough to  proceed as in the construction of the function $U_{\d,\e}^R$ in the proof of Theorem~1, by letting 
$v$ be the solution of the Dirichlet problem $-\D v=F$ on $B_{2R}$ and $v=0$ on $\p B_{2R}$, where 
$$F(x)=\cases{ f  & if $|x|<R$\cr 0 & if $R\le |x|<{4\over3} R$ \cr -\ds{1\over|B_{{5\over3}R}\setminus B_{{4\over3}R}|}\int_{B_R} f & if ${4\over 3}R\le |x|< {5\over3} R$\cr 0& if ${5\over3}R\le |x|\le 2R.$}
$$
In summary, we have that $u^*\le (u_0)^*+\lambda\le (v_0)^*+\l\le v_1^*+\l\le v^*+C+\l$ and this proves our initial claim \eqref{bound} and therefore \eqref{@15}.

Suppose now that $X$ is an r.i. space with the Fatou property, and that $\Wo(\Omega)\subseteq X$. Then $\Lambda_n^0(\Omega)\subseteq X$, continuously, so it is an easy matter to check that when   $u\in L_\expp(\Omega)$ ($n=2$) or $u\in L^{{n\over n-2},\infty}(\Omega)$ ($n\ge 3$), then $u\in X$, by considering the sequence of truncations $u_n=\min\{|u|,n\}$, which belongs to $\Lambda_n^0(\Omega)$, has uniformly bounded norm, and converges monotonically to $|u|$.
\endpf
\smallskip
\noin{\bf Remark.} Estimate \eqref{bound} can be extended to arbitrary functions $u$ in $L_\expp(\Omega)$ ($n=2$) or in $L^{{n\over n-2},\infty}(\Omega)$ ($n\ge 3$) as follows: 
$$u^*(t)\le v^*(t)+C+N_{|B|}^*(t)\;\limsup_{s\to 0} {u^{**}(s)\over N_{|B|}^*(s)},\qquad 0<t\le |\Omega|,\eqdef{bound}$$
for some $v\in \Wo(\Omega)$, some ball $B\subseteq \Omega$ and some constant $C$. This follows from the previous proof, 
since the function $k$ in \eqref k satisfies $k(0)= \limsup_{s\to 0} {u^{**}(s)\over N_{|B|}^*(s)}$.
\bigskip 
\centerline{\bf Acknowledgments}\bigskip
The authors wish to thank Andrea Cianchi and Vincenzo Ferone for pointing out the articles  [Al] and [AF], and Mario Milman for useful discussions. \bigskip

\eject\centerline{\bf References}
\bigskip
\item{[A]} Adams D.R., {\sl A sharp inequality of J. Moser for higher order derivatives}, Ann. of Math. {\bf 128} (1988), 385-398. \smallskip
\item{[AC]} Alberico A., Cianchi A.,
{\sl Optimal summability of solutions to nonlinear elliptic problems},
Nonlinear Anal. {\bf 67 }(2007), 1775-1790. \smallskip
\item{[AF]} Alberico A., Ferone V., {\it Regularity properties of solutions of elliptic equations in $\R^2$ in limit cases}, Atti Accad. Naz. Lincei Cl. Sci. Fis. Mat. Natur. Rend. Lincei (9) Mat. Appl. {\bf 6} (1995),  237–250.\smallskip 
\item{[Al]} Alvino A., {\it A limit case of the Sobolev inequality in Lorentz spaces},
Rend. Accad. Sci. Fis. Mat. Napoli {\bf{ 44}} (1977), 105-112.\smallskip 
\item{[AFT]} Alvino A., Ferone V., Trombetti G., {\sl Estimates for the gradient of solutions of nonlinear elliptic equations with $L^1$ data}, Ann. Mat. Pura Appl {\bf 178} (2000), 129-142.\smallskip  \item{[BM]}
Brezis H., Merle F., {\sl Uniform estimates and blow-up behavior for solutions of $-\D u = V(x)e^u$ in two dimensions}, 
Comm. Part. Diff. Eq. 16 (1991),  1223-1253.\smallskip 
\item{[BS]} Bennett C., Sharpley R., {\sl Interpolation of Operators}, Pure and  Applied  Math., 129, Academic Press, Inc., Boston, MA, 1988.\smallskip
\item{[BSt]} Brezis H., Strauss W.A., {\sl Semi-linear second-order elliptic equations in $L^1$}, 
J. Math. Soc. Japan {\bf25} (1973), 565-590.\smallskip 
\item{[BW]} Brezis H., Wainger S., {\sl A note on limiting cases of Sobolev embeddings}, Comm. Partial Differential Equations {\bf 5} (1980), 773-789.\smallskip
\item{[C]} Cianchi A., {\sl Higher-order Sobolev and Poincar\'e inequalities in Orlicz spaces}, 
Forum Math. {\bf 18} (2006),  745-767.\smallskip 
\item{[CP]} Cwikel M., Pustylnik E., {\sl Sobolev type embeddings in the limiting case}, J. Fourier Anal. Appl. {\bf 4} (1998), 433-446.\smallskip
\item{[CRT]} Cassani D., Ruf B., Tarsi C., {\sl Best constants in a borderline case of second-order Moser Type inequalities}, Ann. Int. H. Poincar\'e Anal. Non Lin\'eaire  {\bf 27} (2010) 73-93.\smallskip
\item{[EKP]} Edmunds D.E., Kerman R.,  Pick L., {\sl Optimal Sobolev imbeddings involving rearrangement-invariant quasinorms}, J. Funct. Anal. {\bf 170} (2000) 307-355.\smallskip

\item{[H]}  Hansson K., {\sl Imbedding theorems of Sobolev type in potential theory}, Math. Scand. {\bf 45} (1979), 77-102.\smallskip
\item{[KL]} Kaminska A., Lee H.J.,
{\it M-ideal properties in Marcinkiewicz spaces},  
Comment. Math. Prace Mat. Tomus specialis in Honorem Juliani Musielak, (2004), 123-144.\smallskip 
\item{[KPS]}  Krein S.G., Petunin J.I., Semenov E.M., {\sl Interpolation of linear operators},  Translations of Mathematical Monographs, 54. American Mathematical Society, Providence, R.I., 1982.\smallskip
\item{[O]} Ornstein D., {\sl A non-equality for differential operators in the $L_1$ norm}, Arch. Ration. Mech. Anal. {\bf 11} (1962), 40-49.\smallskip
\item{[MP]} Milman M.,  Pustylnik E., {\sl On sharp higher order Sobolev embeddings}, Commun. Contemp. Math. {\bf6} (2004),  495-511.\smallskip
\item{[T]} Talenti G., {\sl Elliptic Equations and rearrangements}, Ann. Scuola Norm. Sup. Pisa {\bf 3} (1976), 697-718.
\bigskip

\noin Luigi Fontana \hskip19em Carlo Morpurgo

\noin Dipartimento di Matematica ed Applicazioni \hskip5.5em Department of Mathematics 
 
\noin Universit\'a di Milano-Bicocca\hskip 13em University of Missouri, Columbia

\noin Via Cozzi, 53 \hskip 19.3em Columbia, Missouri 65211

\noin 20125 Milano - Italy\hskip 16.6em USA 
\smallskip\noin luigi.fontana@unimib.it\hskip 15.3em morpurgoc@missouri.edu

\end